\documentclass{ifacconf}

\usepackage{amsmath,amssymb,amsfonts}
\usepackage{graphicx}      
\usepackage{natbib}        
\usepackage{textcomp}
\usepackage{mathtools}
\usepackage{import}
\usepackage{xcolor}
\usepackage{tikz}
\usepackage{caption}
\usepackage{subcaption}
\usepackage{epstopdf}
\usepackage{dsfont}
\usepackage{algorithm}
\usepackage[]{algpseudocode}

\begin{document}
\begin{frontmatter}

\title{Stochastic MPC with Distributionally Robust Chance Constraints}
\thanks[footnoteinfo]{© 2020 the authors. This work has been accepted to IFAC for publication under a Creative Commons Licence CC-BY-NC-ND.}

\author[First]{Christoph Mark}
\author[First]{Steven Liu} 

\address[First]{Institute  of  Control  Systems,  Department  of  Electrical  and  Computer
Engineering, University of Kaiserslautern, 67663
Kaiserslautern, Germany (e-mail: \{mark, sliu\}@eit.uni-kl.de).}

\begin{abstract}
In this paper we discuss distributional robustness in the context of stochastic model predictive control (SMPC) for linear time-invariant systems. We derive a simple approximation of the MPC problem under an additive zero-mean i.i.d. noise with quadratic cost. Due to the lack of distributional information, chance constraints are enforced as distributionally robust (DR) chance constraints, which we opt to unify with the concept of probabilistic reachable sets (PRS). For Wasserstein ambiguity sets, we propose a simple convex optimization problem to compute the DR-PRS based on finitely many disturbance samples. The paper closes with a numerical example of a double integrator system, highlighting the reliability of the DR-PRS w.r.t. the Wasserstein set and performance of the resulting SMPC. 
\end{abstract}

\end{frontmatter}
\newpage
Model Predictive Control (MPC) is an optimization based control strategy that has received a lot attention through academia and industry over the last couple of decades, mainly due to its ability to explicitly take constraints into account \citep{rawlings2009model}. 

In almost every practical application, the states of the system are falsified by disturbances, which result from e.g. measurement noise, process noise or model-plant mismatch. If one can give an a-priori bound on such disturbances, one may use a robust MPC framework \citep{mayne2005robust}, which guarantees robust constraint satisfaction and robust stability. However, these approaches tend to be overly conservative, since the robustness is regarding worst-case disturbances, which may occur to a very low percentile. 

If one has a model of the disturbance, e.g. the probability distribution, then the designer may use a stochastic MPC framework \citep{mesbah2016stochastic}, which aims to find a trade-off between constraint satisfaction and controller performance. By constraint satisfaction we refer to chance constraints \citep{farina2016stochastic}, which are a popular type of constraint relaxation technique, which enforces that the hard constraints have to be satisfied only with a certain probability. In many practical applications these kind of constraints are very beneficial, e.g. comfort constraints in a room temperature control scenario.

Stochastic MPC approaches can furthermore be separated into several sub categories. On the one hand there are \textit{scenario-based approaches} \citep{schildbach2012randomized} \citep{hewing2019scenario}, which, at every time instant, sample sufficiently many disturbance realizations in order to approximate the stochastic control problem. The inherent sampling technique makes the scenario-based methods applicable for systems with arbitrary disturbances. However, due to their heavy computational load these methods are still limited to small scale systems.
On the other hand we have \textit{analytical approximation methods} \citep{farina2013probabilistic} \citep{hewing2018recursively} \citep{mark2019distributed} \citep{mark2020stochastic}, which assume a certain structure in the probability distribution in order to reformulate the stochastic problem based on the moments of the disturbance. These approaches are limited due to the necessity of the exact probability distribution, which, from a practical point of view, is usually not available.

Motivated by recent results on distributionally robust (DR) optimization \citep{esfahani2018data}, we try to seek a middle path in between scenario-based and analytical approximation methods, where we utilize Wasserstein ambiguity sets \citep{kuhn2019wasserstein} in combination with indirect feedback SMPC \citep{hewing2018recursively}.

\subsubsection*{Related work:}
In \citep{van2015distributionally} the authors propose DR receding horizon and infinite horizon controllers for linear discrete time-invariant (DLTI) systems subject to additive disturbances. The DR chance constraints are approximated with conditional value-at-risk (CVaR) constraints and a moment-based ambiguity set. In \citep{darivianakis2017power} the authors propose to use statistical hypothesis testing to construct probabilistic confidence bounds that replace the DR chance constraints as robust constraints.
\subsubsection*{Contribution:}
In this paper we extend the framework of PRS to DR-PRS. This is achieved by rewriting the PRS by means of the union of convex loss functions and computing for each loss function the worst-case VaR, which is equivalent to the DR chance constraint formulation. For Wasserstein ambiguity sets we propose a data driven method to compute DR-PRS based on finite samples using a CVaR optimization problem of the error loss function.

\subsubsection*{Outline:}
The paper is organized as follows. In the first section we introduce the DR problem setting. Section \ref{sec:DR_SMPC} introduces the MPC ingredients and the core of the paper, namely the reformulation of DR chance constraints as DR-PRS. In section \ref{sec:tractable} we present tractable results for worst-case VaR optimization problems whereas Section \ref{sec:examples} closes the paper with two numerical examples on DR-PRS and an MPC implementation of a second order system. For the sake of readability all proofs can be found in the appendix.

\subsubsection*{Notation:}
Given two polytopic sets $\mathbb{A}$ and $\mathbb{B}$, the Pontryagin difference is given as $\mathbb{A} \ominus \mathbb{B} = \{ a \in \mathbb{A} : a+b \in \mathbb{A}, \forall b \in \mathbb{B} \}$. The set of symmetric positive semidefinite (definite) matrices is defined as $\mathbb{S}, \mathbb{S}_{+}$. Positive definite and semidefinite matrices are indicated as $A>0$ and $A\geq0$, respectively. For an event $E$ we define the probability of occurrence as $\mathbb{P}(E)$, whereas $\mathbb{P}(E|D)$ denotes the conditional probability. We define for a random variable $w$ the expected value w.r.t. distribution $\mathbb{P}$ as $\mathbb{E}_{\mathbb{P}}(w)$. Given a measurable loss function $l(w)$ and a risk threshold $\epsilon \in (0, 1)$, the VaR of $l(w)$ at level $\epsilon$ is defined as $\mathbb{P}\text{-VaR}_\epsilon(l(w)) = \text{inf} [ \gamma \in \mathbb{R} | \mathbb{P} \{ \gamma \leq l(w) \} \leq \epsilon ]$. The CVaR of $l(w)$ at level $\epsilon$ is given as $\mathbb{P}\text{-CVaR}_\epsilon(l(w)) = \text{inf} [\tau + \epsilon^{-1}  \mathbb{E}_{\mathbb{P}} \{ \text{max}(0, l(w) - \tau) \}| \tau \in \mathbb{R}]$.


\section{Problem formulation}
\label{sec:problem}
We consider a discrete-time linear time-invariant system
\begin{align}
	x(k+1) = A x(k) + B u(k) + w(k) \label{eq:real_system}
\end{align}
with state $x \in \mathbb{R}^{n_x}$, input $u \in \mathbb{R}^{n_u}$ and independent and identically distributed (i.i.d.) white noise $w \in \mathbb{R}^{n_x}$. The states and inputs are subject to chance constraints
\begin{subequations}
\begin{align}
	&\mathbb{P}(x \in \mathbb{X} | x(0)) \geq p_x \\
	&\mathbb{P}(u \in \mathbb{U} | x(0)) \geq p_u,
\end{align}
\label{eq:chance_constraints}
\end{subequations}
where the constraint sets $\mathbb{X} = \{ x \in \mathbb{R}^{n_x} | H x(k) \leq h \}$ and $\mathbb{U} = \{ u \in \mathbb{R}^{n_u} | L u(k) \leq l \}$ are convex polytopes that contain the origin in their interior with $h \in \mathbb{R}_{> 0}^{n_p}$ and $l \in \mathbb{R}_{> 0}^{n_q}$. The level of chance constraint satisfaction is regulated with $p_x, p_u \in (0, 1)$.

A common assumption in analytical SMPC approaches is that the statistics of the disturbance $w$ are known and so the probability measure $\mathbb{P}$. However, this is restrictive from a practical point of view, where the true distribution is usually unknown, i.e. it must be estimated from historical data. We make the following assumption.
\begin{assum}~
\begin{itemize}
 \item The true probability distribution $\mathbb{P}$ is light-tailed
 \item $\mathbb{P}$ belongs to an ambiguity set $\hat{\mathcal{P}}$ with probability $1 - \beta$, where $\beta \in (0,1)$.
\end{itemize}
 \label{assum:ccu}
\end{assum}
\begin{rem}
An ambiguity set is a collection of plausible variations of the empirically estimated distribution, where the literature distinguishes roughly between discrepancy and moment-based ambiguity sets. We refer the interested reader to the recently published survey paper \citep{rahimian2019distributionally} for a detailed overview. 
\end{rem}
To immunize \eqref{eq:chance_constraints} against distributional ambiguity, we state the following DR chance constraints
\begin{subequations}
\begin{align}
	&\!\inf_{\mathbb{Q} \in \hat{\mathcal{P}}} \mathbb{Q}(x \in \mathbb{X} | x(0)) \geq p_x \label{eq:DR_cc_state} \\
	&\!\inf_{\mathbb{Q} \in \hat{\mathcal{P}}} \mathbb{Q}(u \in \mathbb{U} | x(0)) \geq p_u.
\end{align}
\label{eq:DR_chance_constraints}
\end{subequations}
The control objective is to approximately minimize the infinite horizon cost
\begin{align}
	J_\infty = \lim_{t \rightarrow \infty} \mathbb{E}_{\mathbb{P}} \bigg\{ \frac{1}{t} \sum_{k = 0}^t l(x(k), u(k)) \bigg \}, \label{eq:infinite_cost}
\end{align}
where $l(x(k), u(k))$ denotes a stage cost function. To this end, we split system \eqref{eq:real_system} into a deterministic nominal part $z(k)$ and a stochastic error part $e(k)$, such that $x(k) = z(k) + e(k)$. Following the lines of classic tube-based MPC, we introduce an affine error feedback controller $u(k) = v(k) + K e(k)$, resulting in the decoupled dynamics
\begin{subequations}
\begin{align}
	z(k+1) &= A z(k) + B v(k) \\
	e(k+1) &= (A + BK) e(k) + w(k),  \label{eq:closed_loop_error_system}
\end{align}
\label{eq:nominal_closed_loop}
\end{subequations}
where $K$ is the solution of the infinite horizon linear-quadratic control problem, such that $A_K \coloneqq A+BK$ is Schur stable. In this paper we use the recently proposed indirect feedback SMPC scheme from \citep{hewing2018recursively}, where the main advantage over classical SMPC approaches is the validity of \eqref{eq:closed_loop_error_system} in closed-loop under the MPC control law. For the treatment of DR chance constraints, we extend the concept of PRS to its DR counterpart for the error system \eqref{eq:closed_loop_error_system}.

\section{Distributionally Robust SMPC}
\label{sec:DR_SMPC}
In this section we introduce the controller structure, the MPC objective function and reformulate the DR chance constraints \eqref{eq:DR_chance_constraints} in terms of a DR-PRS.
\subsection{Controller structure}
In the following we introduce the predictive dynamics for the MPC optimization problem, where we denote predicted quantities as e.g. $x(t|k)$ for a $t$-step ahead prediction of the state initialized at time step $k$. The predictive dynamics are
\begin{subequations}
\begin{align}
	&x(t+1|k) = A x(t|k) + B u(t|k) + w(t|k) \\
	&z(t+1|k) = A z(t|k) + B v(t|k) \label{eq:nominal_dynamics} \\
	&e(t+1|k) = (A+BK) e(t|k) + w(t|k), \label{eq:error_dynamics}
\end{align}
\end{subequations}
which are coupled to the closed-loop dynamics \eqref{eq:real_system}, \eqref{eq:nominal_closed_loop} at each time step $k$ with  $x(0|k) = x(k)$, $z(0|k) = z(k)$, $e(0|k) = e(k)$. The closed-loop error $e(k)$ is unaffected by the choice of $z(0|k) = z(k)$ and thus, evolves autonomously according to \eqref{eq:closed_loop_error_system} \citep{hewing2019scenario}. For a quick overview, the following list gives insight in which quantities serve which purpose.
\begin{itemize}
	\item The nominal dynamics \eqref{eq:nominal_dynamics} serve as a prediction model for the MPC
	\item The predictive disturbance $w(t|k)$ and predictive error $e(t|k)$ are used to optimize the MPC cost
	\item The disturbance $w(k)$ and closed-loop error $e(k)$ are used to compute PRS for constraint tightening	
\end{itemize} 
The resulting input for system \eqref{eq:real_system} is 
\begin{align}
	u(k) = v^*(0|k) + K e(0|k),	\label{eq:tube_controller}
\end{align}
where $v^*(0|k)$ denotes the first element of the optimal input trajectory. 
\subsection{Objective}
In order to minimize the infinite horizon cost \eqref{eq:infinite_cost}, we consider a receding horizon strategy and solve the problem over a finite horizon of length $N$, where a terminal cost $V_f(x)$ is used to approximate the infinite horizon tail. The expected cost can then be written as
\begin{align}
	\mathbb{E}_{{\mathbb{Q}}} \{ {J}_k(x, u) \} = \mathbb{E}_{{\mathbb{Q}}} \bigg \{ V_f(x(N|k) + \sum_{t=0}^{N-1} l(x(t|k), u(t|k)) \bigg \}, \label{eq:stage_cost_expected}
\end{align}
where the expectation $\mathbb{E}_{{\mathbb{Q}}}(\cdot)$ is taken w.r.t. the distribution ${\mathbb{Q}}$. A general way to solve optimization problems with expected cost functions akin to \eqref{eq:stage_cost_expected} are scenario approaches, e.g. \citep{hewing2019scenario}, where the stochastic program is approximated with a \textit{sampling-average-approximation} (SAA). The idea is  to collect $M$ data samples $\mathcal{W} = \{ w^{(1)}, \ldots, w^{(M)} \}$ and replace $\mathbb{Q}$ with the empirical distribution $\hat{\mathbb{P}} \coloneqq M^{-1}\sum_{j \in \mathcal{W}} \delta_{w^{(j)}}$, where $\delta_{w^{(j)}}$ is the Dirac delta measure concentrated at $w^{(j)}$, such that \eqref{eq:stage_cost_expected} is approximated with $\mathbb{Q} = \hat{\mathbb{P}}$ as
\begin{align*}
\frac{1}{M} \sum_{i=1}^M \bigg \{ V_f(x^{(i)}(N|k) + \sum_{t=0}^{N-1} l(x^{(i)}(t|k), u^{(i)}(t|k)) \bigg \}.
\end{align*}
Minimization of the above cost yields a control sequence $u^*(\cdot|k)$ that provides an upper bound of the \textit{in-sample performance} for the data set $\mathcal{W}$. However, if we implement $u^*(0|k)$ to the real system \eqref{eq:real_system}, which introduces new (disturbance) samples, then the controller may show a poor \textit{out-of-sample performance} $\mathbb{E}_{{\mathbb{P}}} \{ {J}_k(x, u^*) \}$, which is the real quantity of interest. To highlight the importance between these two performance criteria, we refer the reader to \citep{esfahani2018data}.
\begin{rem}
\label{rem:SAA}
The SAA converges almost surely to the true expectation $\mathbb{E}_{\mathbb{P}}\{ \cdot \}$ as the data size $M$ tends to infinity \citep{esfahani2018data}. However, for small data sizes, the SAA performs poorly in regard to the out-of-sample performance and is the main motivation to consider DR optimization.
\end{rem}
To robustify the MPC optimization problem against distributional uncertainty, we consider an ambiguity set $\hat{\mathcal{P}}$ to reformulate the cost \eqref{eq:stage_cost_expected} as 
\begin{align}
 \!\sup_{\mathbb{Q} \in \hat{\mathcal{P}}} \mathbb{E}_{\mathbb{Q}} \bigg \{ V_f(x(N|k) + \sum_{t=0}^{N-1} l(x(t|k), u(t|k)) \bigg \}, \label{eq:dr_cost}
\end{align}
which serves as an upper bound for the out-of-sample performance with high probability, i.e.
\begin{align}
	{\mathbb{P}} \big \{ \mathbb{E}_{{\mathbb{P}}} \{ {J}_k(x, u) \} \leq \eqref{eq:dr_cost} \big \} \geq 1 - \beta. \label{eq:finite_sample_guarantee}
\end{align}
\subsection*{Quadratic cost with zero-mean i.i.d. noise}
Since the main emphasize of this paper are DR chance constraints, we consider the special case of zero-mean i.i.d. disturbances under a quadratic cost function
\begin{align*}
	l(x,u) = \Vert x \Vert^2_Q + \Vert u \Vert^2_R,
\end{align*}
where $Q \geq 0$, $R>0$ are weighting matrices for the states and inputs. The terminal cost function is given by
\begin{align*}
	V_f(x) = \Vert x \Vert_P^2
\end{align*}
with a positive definite weight $P > 0$ that satisfies the Lyapunov equation  $A_K^\top P A_K + Q + K^\top R K = P$.
Note that the closed-loop error $e(k)$ is driven by the stochastic process \eqref{eq:closed_loop_error_system}, which, under the zero-mean assumption on $w$ and the initialization of $e(0) = 0$, is zero-mean. For the predicted error $e(t|k)$ this does not hold in general, which is due to the choice $e(0|k) = e(k)$. Therefore, the predicted state and input mean are given as $\mu_x = z(t|k) + A_K^t e(k)$ and $\mu_u = v(t|k) + K A_K^t e(k)$, such that \eqref{eq:dr_cost} evaluates to
\begin{subequations}
\begin{align}
	&\!\sup_{\mathbb{Q} \in \hat{\mathcal{P}}} \mathbb{E}_{\mathbb{Q}}\bigg \{  \Vert x(N|k) \Vert_P^2 + \sum_{t=0}^{N-1} \Vert x(t|k) \Vert_Q^2 + \Vert u(t|k) \Vert_R^2 \bigg \} \nonumber\\
	&=  \Vert \mu_x(N|k) \Vert_P^2 + \sum_{t=0}^{N-1} \Vert \mu_x(t|k) \Vert_Q^2 + \Vert \mu_u(t|k) \Vert_R^2 \label{eq:nominal_cost_simplified} \\
	&+ \!\sup_{\mathbb{Q} \in \hat{\mathcal{P}}} \mathbb{E}_{\mathbb{Q}} \bigg\{\Vert e(N|k)\Vert_{P}^2 + \sum_{t=0}^{N-1} \Vert e(t|k)\Vert_{Q^*}^2 \bigg \} \label{eq:cost_last_line}
\end{align}
\label{eq:cost_simplified}
\end{subequations}
where $Q^* = Q + K^\top R K$. It can be seen that \eqref{eq:cost_last_line} does not depend on the optimization variables $z,v$ and hence, is a neglectable constant in the optimization problem, which was also reported in \citep{hewing2018recursively}.
\begin{rem}
Since \eqref{eq:cost_last_line} does not depend on the MPC optimization variables $z,v$, the out-of-sample performance \eqref{eq:dr_cost} cannot be improved through the MPC optimization problem. This changes if the underlying disturbance $w$ is non i.i.d. and non zero mean. However, DR performance could still be achieved by considering a DR tube-controller, see \citep{yang2018wasserstein} for related results.
\end{rem}
\subsection{Distributionally robust chance constraints}
The following section is dedicated to the reformulation of the DR chance constraints \eqref{eq:DR_chance_constraints} in terms of a DR Probabilistic Reachable Set.
\begin{defn}
\label{def:dr_prs}
A set $\mathbb{A}$ is a Distributionally Robust Probabilistic Reachable Set (DR-PRS) of probability level $p$ w.r.t. to the ambiguity set $\hat{\mathcal{P}}$ for system \eqref{eq:closed_loop_error_system} initialized with $e(0) = 0$ if
\begin{align*}
	\mathbb{Q}( e(k) \in \mathbb{A}) \geq p \quad \forall \mathbb{Q} \in \hat{\mathcal{P}} \quad  \forall k \geq 0
\end{align*}
\end{defn}
According to \citep{zymler2013distributionally} the DR-PRS condition $\mathbb{Q}( e(k) \in \mathbb{A}) \geq p \quad \forall \mathbb{Q} \in \hat{\mathcal{P}}$ can similarly be rewritten as
\begin{align}
	\!\inf_{\mathbb{Q} \in \hat{\mathcal{P}}} \mathbb{Q}( e(k) \in \mathbb{A}) \geq p \Leftrightarrow \!\sup_{\mathbb{Q} \in \hat{\mathcal{P}}} \mathbb{Q}(e(k) \notin \mathbb{A}) \leq 1-p. \label{eq:equivalence_inf_sup}
\end{align}
\begin{rem}
A DR-PRS is a probabilistic confidence region, such that for all distributional variations over the ambiguity set $\hat{\mathcal{P}}$, any future (unseen) process realization $e(k)$ lies in $\mathbb{A}$ with a probability of at least $p$, i.e. the worst-case probability of $e(k) \notin \mathbb{A}$ is smaller that $1 - p$. 
\end{rem}
\begin{rem}
For time-varying or non zero-mean errors $e$ it is possible to consider time-varying DR-PRS for constraint tightening, i.e. for each time step $k$ we define a DR-PRS $\mathbb{A}_k$, such that 
\begin{align*}
\!\inf_{\mathbb{Q} \in \hat{\mathcal{P}}} \mathbb{Q}(e(k) \in \mathbb{A}_k) \geq p_x \quad \forall \mathbb{Q} \in \hat{\mathcal{P}}
\end{align*} 
\end{rem}
\begin{rem}
If $\hat{\mathcal{P}}$ is a singleton that contains only the true distribution $\mathbb{P}$, then the ambiguous free PRS, or simply PRS, $\mathbb{S}$ is recovered. 
\end{rem}
The DR-PRS allows for a direct reformulation of the closed-loop chance constraints \eqref{eq:chance_constraints}, which is formalized in the following lemma.
\begin{lem}
\label{lem:chance_constraints}
Consider two polytopic DR-PRS for the states and inputs, i.e. $\mathbb{A}_x = \{ H e(k) \leq \eta_x \}$ and $\mathbb{A}_u = \{ L K e(k) \leq \eta_u \}$. System \eqref{eq:real_system} satisfies the DR chance constraints \eqref{eq:DR_chance_constraints} for any $k \geq 0$ conditioned on $e(0) = x(0) - z(0) = 0$ with a probability of at least $1 - \beta$, if and only if the nominal system \eqref{eq:nominal_dynamics} satisfies the constraints $z(k) \in \mathbb{Z}$ and $v(k) \in \mathbb{V}$ with
\begin{align*}
	\mathbb{Z} = \mathbb{X} \ominus \mathbb{A}_x \\
	\mathbb{V} = \mathbb{U} \ominus \mathbb{A}_u
\end{align*}
\end{lem}
The previous result gives insight in the connection between between DR-PRS and PRS, which is summarized in the following corollary.
\begin{cor}
\label{corol}
Let $\mathbb{S}$ be a PRS and $\mathbb{A}$ a DR-PRS and let $\hat{\mathcal{P}}$ be an ambiguity set. If condition $\mathbb{P} ( {\mathbb{P}} \in \hat{\mathcal{P}}) \geq 1 - \beta$ is satisfied, then ${\mathbb{P}} \big \{ \mathbb{S} \subseteq \mathbb{A} \} \geq 1 - \beta$.
\end{cor}
The main challenge in the DR-PRS synthesis poses the underlying worst-case VaR optimization problem \eqref{eq:wc_var}. The worst-case VaR computation depends highly on the structure of the ambiguity set, e.g. if $\hat{\mathcal{P}}$ is replaced with a moment-based ambiguity set, then the results from \citep{zymler2013distributionally} can be used to express a DR-PRS. A similar result is reported by \citep{darivianakis2017power}, where the authors propose to use statistical hypothesis theory to define a compact DR confidence region to replace a moment-based ambiguity set. In this paper, we opt to solve the worst-case VaR over a Wasserstein ambiguity set, which is presented in Section \ref{sec:wc_var}.
\subsection{Indirect feedback SMPC}
Using the results from the previous sections, we can state the following MPC optimization problem.
\begin{subequations}
\begin{alignat}{2}
&\!\min_{v}  &   &  \quad  \eqref{eq:cost_simplified}  \label{eq:cost_mpc_prob} \\
&\text{s.t.} &   & \quad \mu_x(t+1|k) = z(t+1|k) + e(t+1|k), \\
&            &   & \quad z(t+1|k) = A z(t|k) + B v(t|k),  \\
&            &   & \quad e(t+1|k) = (A+BK)^t e(k),  \\
&            &   & \quad \mu_u(t|k) = v(t|k) + K e(t|k),  \\
& 			 &   & \quad  z(t|k) \in \mathbb{Z}, v(t|k) \in \mathbb{V}, z(N|k) \in \mathbb{Z}_f \\
& 			 &   & \quad  \mu_x(0|k) = x(k), z(0|k) = z(k), e(0|k) = e(k)
\end{alignat}
\label{eq:mpc_prob_tractable}
\end{subequations}
for all $t = 0, \ldots, N-1$. For $k=0$ the nominal state is initialized with $z(0) = x(0)$. The result of the optimization problem is a sequence of nominal control inputs $v(t|k), t = \{0, \ldots, N-1\}$, whereas only the first element $v^*(0|k)$ is implemented with control law \eqref{eq:tube_controller} to the real system \eqref{eq:real_system} and the rest is discarded. Then the time step $k \leftarrow k+1$ increments and problem \eqref{eq:mpc_prob_tractable} is solved repeatedly. For the sake of simplicity we assume a zero terminal constraint, i.e. $\mathbb{Z}_f = \{0\}$. Recursive feasibility of problem \eqref{eq:mpc_prob_tractable} can be established with the standard procedure (shifted optimal cost and terminal controller), see \citep[Theorem 1]{hewing2018recursively} for details.
Satisfaction of the DR chance constraints (with $1-\beta$ confidence) in closed-loop follows from Lemma \ref{lem:chance_constraints} and \citep[Lemma 1]{hewing2018recursively}, which states that the predicted error has the same distribution as the closed-loop error.
\section{Tractable reformulations}
\label{sec:tractable}
The general ambiguity set $\hat{\mathcal{P}}$ is replaced with a Wasserstein ambiguity set, which is defined as a statistical ball of radius $\theta$ centered at the empirical distribution $\hat{\mathbb{P}}$, i.e.
\begin{align}
	\mathbb{B}_\theta(\hat{\mathbb{P}}) = \{ \mathbb{Q} \in \mathcal{P}(\mathbb{R})\: |\: d_p(\hat{\mathbb{P}}_M, \mathbb{Q}) \leq \theta \}, \label{eq:ambiguity_set}
\end{align}
where $\mathcal{P}(\mathbb{R})$ is a set of Borel probability measures on $\mathbb{R}$, such that the expected value $\mathbb{E}_{\mathbb{Q}}(\Vert w \Vert_p)$ under any norm $\Vert \cdot \Vert_p$ is finite. The measure $ d_p(\mathbb{Q}, \mathbb{P})$ denotes the $p$-Wasserstein distance.
\begin{defn}
For any $p \in [1, \infty)$, the $p$-Wasserstein distance between two arbitrary distributions $\mathbb{Q}$ and $\mathbb{P}$ supported on $\mathbb{R}$ is defined as
\begin{align*}
	d_p(\mathbb{Q}, \mathbb{P}) \coloneqq \!\inf_{\Xi}  \bigg( \int_{\mathbb{R}^2} \Vert \xi - \xi^\prime \Vert_p \Xi(d \xi, d \xi^\prime) \bigg)^{\frac{1}{p}} ,
\end{align*}
where $\Xi$ is the set of all joint probability distribution of $\xi$ and $\xi^\prime$ with marginals $\mathbb{Q}$ and $\mathbb{P}$.
\end{defn}
A Wasserstein radius $\theta$ for \eqref{eq:ambiguity_set} that satisfies Assumption \ref{assum:ccu} can is formalized in the following theorem.
\begin{thm}
Let Assumption \ref{assum:ccu} hold with $\beta \in (0,1)$ and let $\hat{\mathcal{P}} = \mathbb{B}_\theta(\hat{\mathbb{P}})$. Then, for a given data size $M$, there exists an $\theta = \theta_M(\beta) > 0$, such that \eqref{eq:finite_sample_guarantee} holds.
\end{thm}
The proof of the above theorem follows from \citep[Thm. 3.5]{esfahani2018data} and more generally by \citep[Thm. 2]{fournier2015rate}. The previous result highlights the main advantage of the Wasserstein ambiguity set compared to moment-based ambiguity sets, since it connects the data size $M$ to the size of the uncertainty set while maintaining the probabilistic guarantee $1-\beta$. Furthermore, the Wasserstein radius $\theta_M(\beta)$ is a monotonic decreasing function in the data size $M$, i.e. for $M \rightarrow \infty$, the Wasserstein ambiguity set reduces to the singleton that contains only the true distribution $\mathbb{P}$. 
\subsection{Worst-case Value-at-Risk}
\label{sec:wc_var}
\begin{figure}[htbp]
\centering
	  \includegraphics[width=0.7\linewidth]{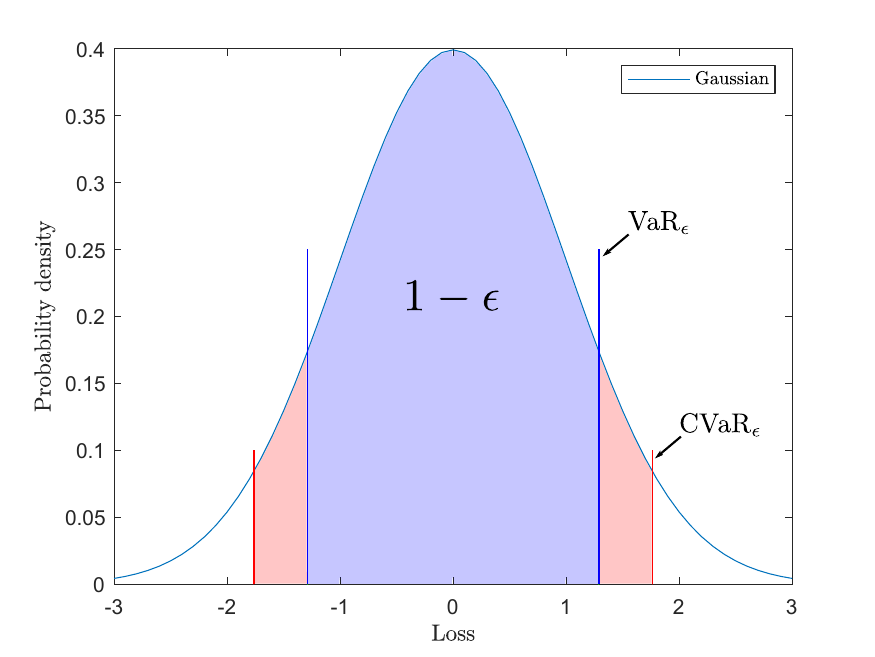}
	  \caption{Box-shaped PRS with VaR and CVaR for a standard Gaussian distribution at probability level $p = 0.8$. } 
	  \label{fig:cvar}
\end{figure}

In this section we introduce a method to solve worst-case VaR problems that appear in the DR-PRS synthesis. We consider a set $\mathcal{E} = \{e^{(1)}, \ldots, e^{(M)} \}$ that contains $M$ i.i.d. samples of the closed-loop error \eqref{eq:closed_loop_error_system} and use the marginal distribution in each direction of the error $e$ to define a symmetric box-shaped DR-PRS of the form
\begin{align*}
	\mathbb{A}_B \coloneqq \{ e \in \mathbb{R}^{n_x} \: | \: |e_i|  \leq \eta_i \quad \forall i = 1, \ldots, n_x\},
\end{align*}
where $e_i$ denotes the $i$-th element of $e$. Using techniques from \citep{esfahani2018data}, the worst-case VaR $\eta_i$ of the $i$-th marginal distribution at probability level $\epsilon_i$ can be computed as a byproduct of the following worst-case CVaR optimization problem
\begin{subequations}
\begin{alignat}{2}
\text{CVaR}_{\epsilon_x}(|e_i|) = \: \ &\!\min_{\lambda \geq 0}  &   & \:\: \tilde{\eta}_i + \frac{1}{\epsilon_i} \bigg ( \theta \lambda + \frac{1}{M} \sum_{j=1}^M \alpha_j \bigg )\\
&\text{s.t.} &   & \quad \theta \lambda + \frac{1}{M} \sum_{j=1}^M \alpha_j \leq \epsilon_i \label{eq:Wasserstein_dependence}\\
&            &   & \quad \alpha_j \geq ( |e_i^{(j)}| - \tilde{\eta}_i)^+   \hspace{0.35em} \forall j \in \mathcal{E} \label{eq:non_convex}
\end{alignat}
\label{eq:tractable_wasserstein}
\end{subequations}
where $(|e_i^{(j)}| - \tilde{\eta}_i)^+ = \text{max}(0, |e_i^{(j)}| - \tilde{\eta}_i)$. Note that for $\theta = 0$, the optimization problem is equal to a SAA, i.e. the worst-case CVaR and worst-case VaR reduce to the CVaR and VaR based on the empirical distribution $\hat{\mathbb{P}}$. Furthermore, the Wasserstein penalty $\lambda$ can be increased arbitrarily without having an effect on $\eta_i$.
\begin{rem}
The CVaR is a measure for the average loss above the $(1-\epsilon_i)$-th quantile of the loss function, i.e. $\text{VaR}_{\epsilon}(|e_i|) \leq \text{CVaR}_{\epsilon}(|e_i|)$,
as illustrated in Figure \ref{fig:cvar}.
\end{rem}
\section{Numerical example}
\label{sec:examples}
This section is dedicated to a numerical example, highlighting the convergence and chance constraint satisfaction properties of the proposed approach. To compare our results, we utilize the double integrator example from \citep{hewing2018recursively} with dynamics
\begin{align*}
	x(k+1) = \begin{bmatrix}
		1 & 1 \\
		0 & 1
	\end{bmatrix} x(k) + \begin{bmatrix}
		0.5 \\
		1
	\end{bmatrix} u(k) + w(k),
\end{align*}
where $w$ is a normally distributed zero-mean i.i.d. disturbance with true distribution
\begin{align*}
	w \sim \mathcal{N}(0, \Sigma_w), \quad \Sigma_w = \begin{bmatrix}
	0.25 & 0.5\\
	0.5 & 1
	\end{bmatrix}.
\end{align*}
We impose a chance constraint on the second state, i.e. $\mathbb{P}(| x_2(k) | \leq 1.2 ) \geq 0.8$. Furthermore, we chose the stabilizing controller $K = [-0.2, -0.6]$. The cost weighting matrices are given as $Q = I$, $R = 1$ and the prediction horizon is $N = 30$. For the sake of simplicity we use a zero terminal constraint, i.e. $\mathbb{Z}_f = \{ 0 \}$.
\subsection{Ambiguity set}
\begin{figure}[h]
\centering
	  \includegraphics[width=0.77\linewidth]{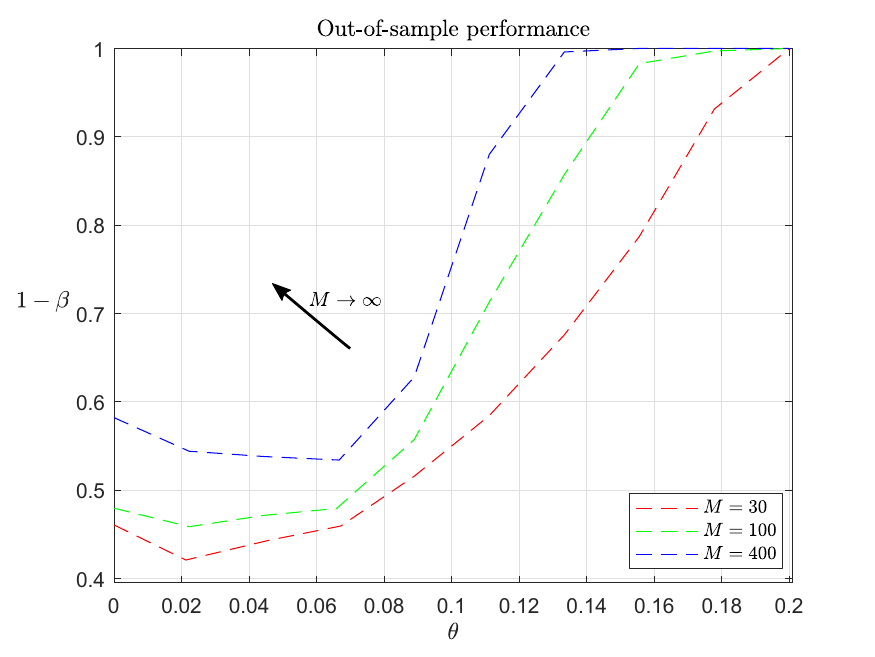}
	  \caption{DR-PRS out-of-sample performance for many Wasserstein radii and data sizes $M = \{30, 100, 400 \}$.} 
	  \label{fig:theta}
\end{figure}

\begin{figure}[h]
\centering
	  \hspace{-0.82em}\includegraphics[width=0.7\linewidth]{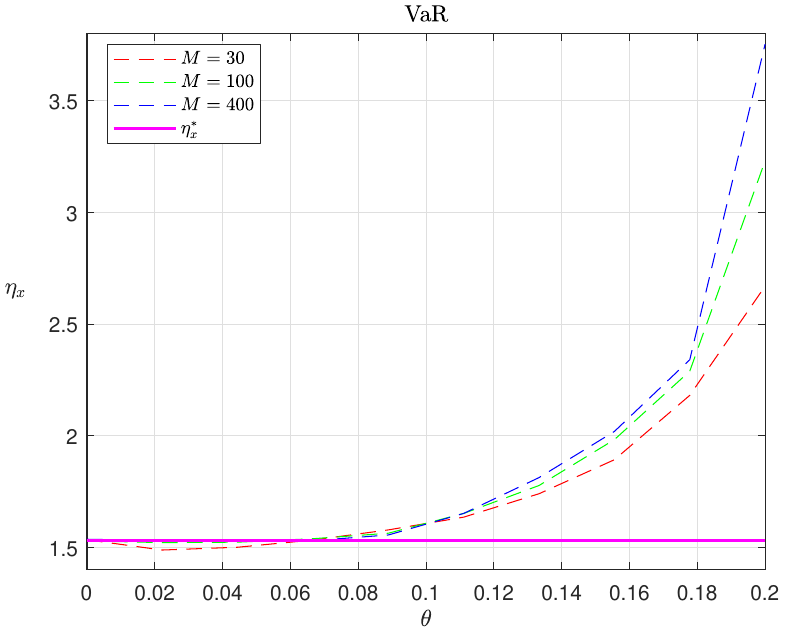}
	  \caption{Quantiles for the data sizes $M = \{30, 100, 400 \}$. The magenta line depicts the true VaR $\eta_x^*$ = 1.53.} 
	  \label{fig:quantiles}
\end{figure}

In the first experiment we want to exemplify the impact of the Wasserstein radius on the reliability of the DR-PRS. First, we compute the true error distribution $\mathbb{P} = \mathcal{N}(0,\Sigma_e)$, where $\Sigma_e$ solves the Lyapunov equation $\Sigma_e = A_K \Sigma_e A_K^\top + \Sigma_w$. Thus, the true PRS is given as
\begin{align*}
	\mathbb{S} = \bigg \{ e_2 \in \mathbb{R} \: | \: |e_2| \leq \sqrt{ [\Sigma_e]_{2,2} \mathcal{X}^2_{1}(0.8)} = 1.53 \bigg\},
\end{align*} 
where $\mathcal{X}^2_{1}(0.8)$ is the chi-squared distribution at probability level $0.8$ with $1$ degree of freedom.
Then we samples $N_v = 10^5$ error states and store them as validation data. Now we consider three different data sizes $M = \{30, 100, 400 \}$, where for each data size, $N_{mc} = 1000$ DR-PRS are computed via \eqref{eq:tractable_wasserstein} with $\epsilon_{x} = 0.2$. 
For a given $\theta$ the reliability of the $i$-th DR-PRS is given by
\begin{align*}
	\tilde{p}_i = \frac{1}{N_v} \sum_{j=1}^{N_v} \mathds{1} \big \{ | e_2^{(j)} | \leq \eta_x(\theta) \big \} \quad \forall i \in \{1, \ldots, N_{mc} \}
\end{align*}
where $\mathds{1}\{\cdot\}$ is the indicator function. We made explicit that $\eta_x(\theta)$ depends on $\theta$. The reliability of the Wasserstein radius $\theta$ w.r.t. the validation data is given by
\begin{align*}
	r(\theta) = \frac{1}{N_{mc}} \sum_{i=1}^{N_{mc}} \mathds{1}  \{ \tilde{p}_i \geq 0.8 \}.
\end{align*}
Using the former empirical procedure, we compute the reliability for many Wasserstein radii in the interval $0 \leq \theta \leq 0.2$ and summarized the results in Figure \ref{fig:theta}.
In Figure \ref{fig:quantiles} we plot the corresponding quantiles over the Wasserstein radius.
In Figure \ref{fig:theta} it can be observed that the reliability of the VaR estimate depends on the data size $M$ and the Wasserstein radius $\theta$. First we consider the ambiguous free case where $\theta = 0$, which equals a sampling average approximation with the training data set. For small sample sizes the variance of $\eta(\theta)$ is large, which can be seen in Figure \ref{fig:quantiles} for $\theta=0$, where the true quantile is approximated empirically via $N_{mc} = 1000$ Monte-carlo simulations, but the reliability of that estimate is merely $46 \%$. For increasing data sizes $M$, the variance of ${\eta}_x(\theta)$ is reduced, e.g. for $M = 400$ the reliability of the same estimate is  $58 \%$. This verifies the statement of Remark \ref{rem:SAA} that the SAA performs poorly for small data sizes.

 For $\theta > 0$ the reliability of the estimate can be controlled, e.g. for $M=30$ a radius of $\theta=0.2$ is mandatory to guarantee a $99\%$ confidence, whereas for $M = 400$ a radius of $\theta = 0.131$ is sufficient. Thus, it is evident that a $1-\beta$ confidence can be maintained while reducing the Wasserstein radius as the data size $M$ increases. In practice, the Wasserstein radius should initially be set to a large value and successively be reduced by observing the number of constraint violations. 
\subsection{SMPC results}
\begin{figure}[h]
\centering
	  \includegraphics[width=0.85\linewidth]{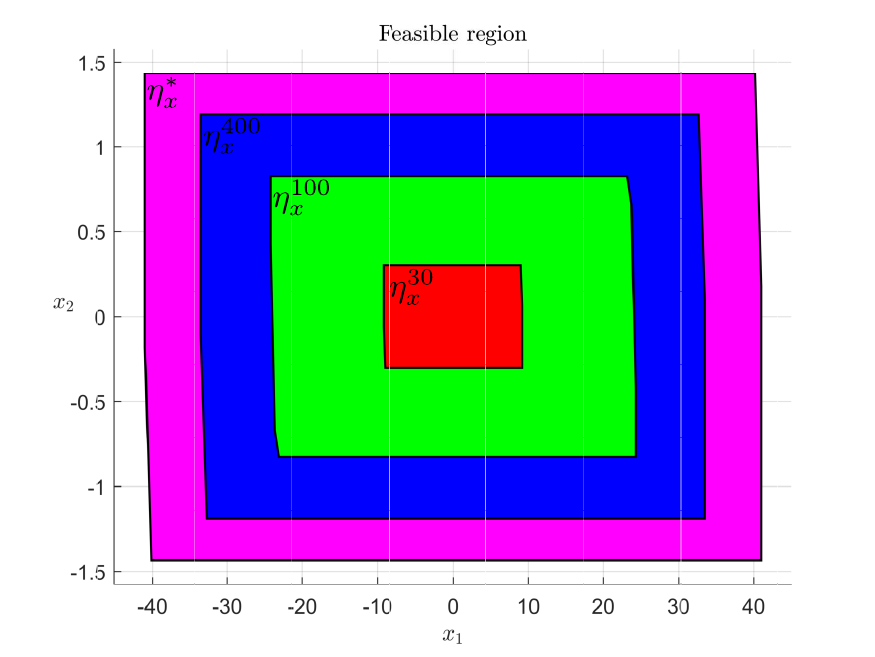}
	  \caption{Feasible regions for different DR-PRS.} 
	  \label{fig:feasible}
\end{figure}

Now we implement three DR-PRS with an empirical reliability of $1 - \beta = 0.99$ with $\eta_{x}^{30} = 2.67$, $\eta_x^{100} = 2.14$, $\eta_x^{400} = 1.8$ and run for each set $N_{s} = 1000$ Monte-carlo simulations of $100$ closed-loop time steps. The initial state is given by $x(0) = [10, 0]^\top$. We compare the average closed-loop cost (Av{[}J{]}) and the number of empirical constraints violations ($\#C_{\text{vio}}$) with the indirect feedback SMPC based on full distributional information with $\eta_x^* = 1.53$.
\begin{center}
  \captionof{table}{Performance for different PRS parameterizations \label{table1}}
\centering
\begin{tabular}{c|cccc}
 PRS       & $\eta_x^{30}$& $\eta_x^{100}$ & $\eta_x^{400}$ & $\eta_x^*$ \\ \hline
Av{[}J{]}    & -         &  9.65        &   8.37       &   7.80   \\ 
$\#C_{\text{vio}}$ & -     &  1091        &   1169       &   1229
\end{tabular}
\end{center}
For $\eta_x^{30}$ the SMPC is infeasible for the initial state, which can be seen in Figure \ref{fig:feasible}, where we compare the feasible regions for each DR-PRS and true PRS.
From Table \ref{table1} it is evident that for increasing data, the average MPC performance will converge to the true cost of around $7.80$. The same argument holds true for the number of constraint violations, which are all above the required level of $p_x = 0.8$.
\section{Conclusion}
In this paper, we presented a SMPC based on indirect feedback with distributionally robust probabilistic reachable sets. We first robustified the MPC cost against distributional ambiguity and observed that for zero-mean i.i.d. disturbances the out-of-sample performance of the MPC cannot be improved through the MPC decision variables. Afterwards, the classical chance constraints got  robustified against distributional ambiguity by means of a DR-PRS. The paper closed with two numerical examples, highlighting the out-of-sample performance of the DR-PRS and the closed-loop performance of the resulting SMPC.

\bibliography{ifacconf}             

\appendix
\section{Proof of Lemma \ref{lem:chance_constraints}}
\begin{pf}
Consider the state constraints \eqref{eq:DR_chance_constraints} in polytopic form $H x(k) \leq h$ conditioned on $x(0)$, which can equivalently be written as
\begin{multline*}
 \!\inf_{\mathbb{Q} \in \hat{\mathcal{P}}} \mathbb{Q}\bigg(\bigcap_{i=1}^{n_p} \big \{ H_i x(k) \leq h_i \big \} \bigg) \geq p_x \\
 \stackrel{\eqref{eq:equivalence_inf_sup}}{\Longleftrightarrow} \!\sup_{\mathbb{Q} \in \hat{\mathcal{P}}} \mathbb{Q}\bigg(\bigcup_{i=1}^{n_p} \big \{ H_i x(k) \geq h_i \big \} \bigg) \leq 1 - p_x,
\end{multline*}
where $H_i$ denotes the $i$-th row of $H$ and $h_i$ the $i$-th element of $h$. By substitution of $x = z + e$ and application of the union bound we obtain
\begin{align*}
 &\!\sup_{\mathbb{Q} \in \hat{\mathcal{P}}} \mathbb{Q}\bigg(\bigcup_{i=1}^{n_p} \big \{ H_i x(k) \geq h_i \big \} \bigg)\\
&\leq \sum_{i=1}^{n_p} \!\sup_{\mathbb{Q} \in \hat{\mathcal{P}}} \mathbb{Q}\{ H_i z(k) + H_i e(k) \geq h_i \} \leq \sum_{i=1}^{n_p} \epsilon_{x,i}.
\end{align*}
The individual violation probabilities $\epsilon_{x,i} \in (0,1)$ satisfy the condition $\sum_{i=1}^{n_p} \epsilon_{x,i} \leq \epsilon_{x} = 1 - p_x$. Furthermore,
\begin{align}
& \!\sup_{\mathbb{Q} \in \hat{\mathcal{P}}} \mathbb{Q} \big \{H_i z(k) + H_i e(k) - h_i \geq 0 \big \} \leq \epsilon_{x,i} \nonumber \\
	\Longleftrightarrow & \!\sup_{\mathbb{Q} \in \hat{\mathcal{P}}} \mathbb{Q}\text{-VaR}_{\epsilon_{x,i}} \big \{ H_i z(k) + H_i e(k) - h_i \big \} \leq 0 \nonumber \\
	\Longleftrightarrow & H_i z(k) \leq h_i - \underbrace{\!\sup_{\mathbb{Q} \in \hat{\mathcal{P}}} \mathbb{Q}\text{-VaR}_{\epsilon_{x,i}} \big \{ H_i e(k) \big \}}_{ \coloneqq \eta_{x,i}} \label{eq:wc_var}
\end{align}
where the first equivalence follows from the connection between the worst-case VaR and DR chance constraints \citep{zymler2013distributionally} and the second equivalence is due to the translational invariance of the VaR. The worst-case VaR in \eqref{eq:wc_var} can be computed as
\begin{subequations}
\begin{alignat}{2}
\eta_{x,i} = \: \:&\!\min_{\tilde{\eta} \in \mathbb{R}} &   &  \quad \tilde{\eta} \\
&\text{s.t.} &   & \quad \!\sup_{\mathbb{Q} \in \hat{\mathcal{P}}} \mathbb{Q}( H_i e \geq \tilde{\eta}) \leq \epsilon_{x,i}. 
\end{alignat}
\label{eq:DR-PRS_opti}
\end{subequations}
By definition, the VaR $\eta_{x,i}$ coincides with the $(1-\epsilon_{x,i})$-th quantile of $H_i e(k)$. The DR-PRS is obtained by defining $\mathbb{A}_x = \{ H e \leq \eta_x \}$.
To summarize, we established that
\begin{align*}
	p_x \leq &\!\inf_{\mathbb{Q} \in \hat{\mathcal{P}}} \mathbb{Q}( e(k) \in \mathbb{A}_x) = \!\inf_{\mathbb{Q} \in \hat{\mathcal{P}}} \mathbb{Q}( x(k) - z(k) \in \mathbb{A}_x) \\
	= & \!\inf_{\mathbb{Q} \in \hat{\mathcal{P}}} \mathbb{Q}( x(k) \in \{z(k)\} \oplus \mathbb{A}_x) \leq \!\inf_{\mathbb{Q} \in \hat{\mathcal{P}}} \mathbb{Q}( x(k) \in \mathbb{X}) 
\end{align*}
where the last inequality follows from $z(k) \in \mathbb{Z}$ and
\begin{align*}
(\underbrace{\mathbb{X} \ominus \mathbb{A}_x}_{\mathbb{Z}}) \oplus \mathbb{A}_x \subseteq \mathbb{X}.
\end{align*}
Recall that the ambiguity set $\hat{\mathcal{P}}$ itself is a random object, such that the true distribution $\mathbb{P}$ belongs only with $1-\beta$ confidence to $\hat{\mathcal{P}}$ (see Assumption \ref{assum:ccu}). Thus, it immediately follows that 
\begin{align*}
{\mathbb{P}} \bigg \{ \!\inf_{\mathbb{Q} \in \hat{\mathcal{P}}} \mathbb{Q}( x(k) \in \mathbb{X}) \geq p_x \bigg \} \geq 1 - \beta.
\end{align*}
This concludes the proof for the DR state chance constraints \eqref{eq:DR_cc_state}. The DR input chance constraints are similarly derived, starting from $L u(k) \leq l$ and separation $L v(k) + L K e(k) \leq l$. \qed
\end{pf}
\section{Proof of Corollary \ref{corol}}
\begin{pf}
Assume that $\mathbb{P} ( {\mathbb{P}} \in \hat{\mathcal{P}}) \geq 1 - \beta = \tilde{\beta}$ is satisfied. 
From the proof of Lemma \ref{lem:chance_constraints} we already know that the DR-PRS $\mathbb{A}$ is the result of a worst-case VaR optimization problem for each halfspace constraint. Now, for each $i \in \{1, \ldots, n_p\}$, we have
\begin{align*}
\mathbb{P} \bigg \{ \mathbb{P}\text{-VaR}_{\epsilon_{x,i}} \{ H_i e(k)\} \leq \!\sup_{\mathbb{Q} \in \hat{\mathcal{P}}} \mathbb{Q}\text{-VaR}_{\epsilon_{x,i}} \big \{ H_i e(k) \} \bigg \} \geq \tilde{\beta},
\end{align*}
which is a direct consequence of the probabilistic guarantee of the ambiguity set.
The claim follows by taking the union bound and changing to set notation. \qed
\end{pf}
\end{document}